\documentclass[11pt,notitlepage,twoside]{article}
\usepackage{latexsym}
\usepackage{amsmath}
\usepackage{amsfonts}
\usepackage{amssymb}

\newcommand{\e }{\varepsilon }

\newcommand{\be}{\begin{equation}}
\newcommand{\ee}{\end{equation}}
\newcommand{\ba}{\begin{align*}}
\newcommand{\ea}{\end{align*}}

\newtheorem{theorem}{Theorem}[section]

\newtheorem{corollary}{Corollary}[section]

\author{Zakaria Boucheche\\
{\footnotesize  Laboratory of Applied Mathematics and Harmonic
Analysis (LR11-ES52)}\\{\footnotesize  (Faculty of Sciences of Gabes.
Gabes University)}}
\title { \Large \textbf{Optimal result involving the Green's function}}
\begin{document}
\date{ }
\maketitle
\begin{abstract} We investigate a borderline between existence and non-existence of positive solution for a nonlinear elliptic equation involving a critical Sobolev exponent in three-dimensional ball. The method is relied on a suitable choice of the functions used to test two natural ingredients for the associated variational problem.
\end{abstract}
\noindent{ \footnotesize {{\it \\\\2020 Mathematics Subject
Classification.}\,\,\, 35J60, 35B33, 35B09, 35J08, 35A01,  49K15.}
}
\section{Introduction and main results}
\def\theequation{1.\arabic{equation}}\makeatother
\setcounter{equation}{0}
 In this paper we study the following nonlinear elliptic partial differential
equation with zero Dirichlet boundary condition
\begin{equation}\label{problem}
\begin{aligned}
-\Delta& u=K(x)u^{5}+\mu u \quad \mbox{on}\,\,B,\\{}& u>0 \quad
\mbox{on} \,\,B,\quad u=0 \quad \mbox{on} \,\,\partial B,
\end{aligned}
\end{equation}where $B:=B(o,\,1)$ is the unit ball of $\mathbb{R}^3$ which is equipped with its Euclidean norm defined, for any $x=(x_1,\,x_2,\,x_3)\in \mathbb{R}^3,$ by $|x|^2=x_1^2+x_2^2+x_3^2,\,\,K(x)$ is a continuous function on $\bar{B}$ and $0<
\mu< \pi^2.$\\

One motivation for the study of this problem comes from \cite[Question 7]{Br1} where the author has given arise an interesting question dealing with the role of the Green's function in establishing an optimal result related to a certain threshold. Influenced by the pioneering paper \cite{Bo1}, we formulate the following natural question that fulfill our present context:\\\\
{\bf Question 1.} How can we demonstrate the necessity of an existence criterion that seems to be optimal?\\\\
In \cite{Bo1} the author of the present paper has implicitly revealed through the proof of \cite[Theorem 1.2]{Bo1} that the regular part of the Green's function of the Laplacian plays a decisive role in confirming the optimality of the condition \cite[(1.5)]{Bo1}. Continuing this positive procedure, we will try to address the question mentioned above through the problem \eqref{problem}. To this end, let us state the following assumptions: Assume that\\\\
$\mathbf{(K_{\eta}^{d_0})}$ $K(x)=f(|x|)$ on $B$ with $f(t)$ is a non-negative and non-increasing $C^1$-function satisfying
\begin{equation}\label{local}f(t)=K(o)+\eta f_1(t),\quad \forall\,\,0\leq t\leq d_0,\end{equation}
where $K(o)> 0,\,\eta \geq 0$ and $1\geq d_0> 0$ are constants, and $f_1(t)$ is independent of $\eta$ with $f_1(0)=0$ and $f_1^{'}(0)\neq 0.$\\\\
In this case, we will refer to the problem \eqref{problem}
as $\mathbf{(BN)}_{\eta}^{d_0}.$\\\\
$\mathbf{(K)}$ $\limsup\limits_{t\rightarrow
0}\frac{f_1'(t)-f_1'(0)}{t}< +\infty.$\\
Because of the local behaviour of $K(x)$ given in $\mathbf{(K_{\eta}^{d_0})},$ the regular part $g_\mu^{B}(x,\,y)$ of the Green's function of the operator $( -\Delta-\mu )$ on $B,$ with Dirichlet boundary condition, will play a decisive role in establishing an existence criterion. We recall that we have from \cite[Remark 2]{Br1}\begin{equation}\label{z00} g_\mu^B(o,\,o)=-\frac{\sqrt{\mu}\cos(\sqrt{\mu})}{4\pi\sin(\sqrt{\mu})}.\end{equation}Our main result is the following:
\begin{theorem}\label{optimal}{\it\,\,Let $0<
\mu< \pi^2.$ Assume that $K(x)$ satisfies the
assumptions $\mathbf{(K)}$ and $\mathbf{(K_{\eta}^1)}.$ Then there exists a constant $\bar{\eta}> 0$
depending on $f_1(t)$ and $K(o)$ such that if $\,\,0\leq \eta\leq
\bar{\eta},$ then the problem $\mathbf{(BN)}_{\eta}^1$ admits a solution
if and only
if\begin{equation}\label{brezis}-\frac{f_1^{'}(0)\eta}{16\pi K(o)}< g_\mu^{B}(o,\,o).\end{equation}}\end{theorem}The
proof of the sufficiency is obtained by the use of a suitable functional for the asssociated variational problem. And then, in order to apply \cite[Theorem 1.1]{Bo2}
we need to test this functional by a function obtained as a projection on the Sobolev space $H_0^1(B)$ by means of the operator $(-\Delta-\mu).$ As well, the necessity of the condition \eqref{brezis} is established through a meticulous test of the Pohozaev-type identity \cite[(2.18)]{Bo1} by a function that involves the quantity $g_\mu^{B}(o,\,o).$\\\\
{\bf Remark 1.1.} By considering a translation and dilatation in the space and using the fact that $g_{\mu/r^2}^{B(a,\,r)}(a,\,a)=(1/r)g_\mu^{B}(o,\,o),$ Theorem \ref{optimal}, with the corresponding suitable changes, remains valid for a general ball $B(a,\,r)$ of center $a$ and radius $r$ in $\mathbb{R}^3.$\\\\
{\bf Remark 1.2.} I think that the condition $\mathbf{(K)}$ is essential. Indeed, a deep analysis as that given to prove \cite[Proposition 1.1]{Bo1} has the advantage to give arise an existence criterion when the conditions $\mathbf{(K)}$ and \eqref{brezis} are not satisfied.\\\\
{\bf Remark 1.3.} The test function used in the proof of \cite[Lemma 1.3]{BrNir} enables us to get the existence criterion$$-\frac{f_1^{'}(0)\eta}{2K(o)}+\frac{\pi^2}{4}< \mu.$$Regarding the condition \eqref{brezis}, there are some values of $\mu$ for the existence of solution for the problem $\mathbf{(BN)_{\eta}^1}$ which are missed. Indeed, an easy calculation shows that, for $\mu=\bigl(-f_1^{'}(0)\eta/2K(o)\bigr)+\pi^2/4,$ we have$$-\frac{f_1^{'}(0)\eta}{16\pi K(o)}< g_{\mu}^{B}(o,\,o)\quad\quad \mathrm{for\,\,}\eta> 0\,\,\mathrm{small\,\,enough}.$$

As a corollary, we state the following curious example where there is no influence of the function $f_1(t)$ on the range of validity of $\eta.$
\begin{corollary}\label{example}{\it\,\,Let $\,\,0<
\mu< \pi^2$ and $K(x)=K(o)-\eta \bigl(|x|+|x|^2\bigr)$ be a non-negative function on $B.$ Then the problem $\mathbf{(BN)}_{\eta}^1$ admits a solution
if and only
if\begin{equation}\label{e.brezis}\frac{\eta}{16\pi K(o)}< g_\mu^{B}(o,\,o).\end{equation}}\end{corollary}
The proof is an immediate consequence of the arguments used in the proof of Theorem \ref{optimal}.\\Regarding this result, naturally one can ask the following:\\\\
{\bf Question 2.} Let $0< \mu < \pi^2$ and $K(x)$ be a function satisfying $\mathbf{(K_{\eta}^1)}$ and $\mathbf{(K)}.$ Is it true that, for any such $\eta,$ the condition \eqref{brezis} is necessary to obtain a solution for the problem $\mathbf{(BN)_{\eta}^1}?$\\

We finish this section by dealing with the constant $d_0$ given in the condition $\mathbf{(K_{\eta}^{d_0})}.$ In fact, the value of $\bar{\eta}$ for which Theorem \ref{optimal} remains valid has a problem of dependence not only on the function $f_1(t)$ but also on the interval of validity of the expression \eqref{local}. Namely we have the following result:\begin{theorem}\label{optimal2}{\it\,\,Let $0<
\mu< \pi^2.$ Assume that $K(x)$ satisfies the
assumptions $\mathbf{(K)}$ and $\mathbf{(K_{\eta}^{d_0})}.$ Then there exists a constant $\bar{\eta}> 0$
depending on $f_1(t),\,d_0$ and $K(o)$ such that if $\,\,0\leq \eta\leq
\bar{\eta},$ then the problem $\mathbf{(BN)}_{\eta}^{d_0}$ admits a solution
if and only
if\begin{equation}\label{brezis2}-\frac{f_1^{'}(0)\eta}{16\pi K(o)}< g_\mu^{B}(o,\,o).\end{equation}}\end{theorem}The proof is similar, up to minor modifications, to that given for Theorem \ref{optimal}. So we omit it.
\section{Proof of Theorem \ref{optimal} }\label{sec2}
\def\theequation{2.\arabic{equation}}\makeatother
\setcounter{equation}{0}
{\it Sufficiency of the condition \eqref{brezis}:} We define on $$M_5:=\{u\in H_0^1(B):\,\bigl(\int_B|\nabla u|^2\bigr)^{\frac{1}{2}}>6^{-1}\,\,\mathrm{and}\,\,\int_BK(x)u^6> K(o)(5S)^{-3}\bigl(\int_B|\nabla u|^2-\mu \int_Bu^2\bigr)^3\}$$ the functional $J(u)$ by
\begin{equation}\label{bahri}J(u)=\frac{\int_{B}|\nabla
u|^2- \mu \int_{B} u^2} {\Bigl( \int_{B}K(x)u^6 \Bigr)^{\frac{1}{3}}}.
\end{equation}In order to apply \cite[Theorem 1.1]{Bo2}, we have to look for a test function $u$ in $M_5$ such that \begin{equation}\label{lions}J(u)<
\frac{1}{[K(o)]^{\frac{1}{3}}}S,
\end{equation} where $S$ denotes the best Sobolev constant of the embedding $H_0^1(B)$ in $L^6(B).$ To this end, we argue as in the proof of \cite[Theorem 7]{Br1}: Let, for any $\e> 0,\,\,\phi_{\e}$ be the solution of the problem
\begin{equation*}
\begin{aligned}
-\Delta& \phi_{\e}-\mu \phi_{\e}=-\Delta U_{\e} \quad \mbox{on}\,\,B,\\{}& \quad \phi_{\e}=0 \quad \mbox{on} \,\,\partial B,
\end{aligned}
\end{equation*}
where $U_{\e}(x)$ is defined by$$U_{\e}(x)=\frac{\sqrt{\e}}{\bigl(\e^2+|x|^2\bigr)^{\frac{1}{2}}}$$
as a solution of $-\Delta U_{\e}=3U_{\e}^5$ on $\mathbb{R}^3.$ We claim that, for $\e >0$ small enough,\begin{equation}\label{boucheche}\phi_{\e}\in M_5\quad\mathrm{and}\quad J(\phi_{\e})=\frac{S}{[K(o)]^{\frac{1}{3}}}\Bigl(1+\bigl[\frac{-f_1^{'}(0)\eta}{16\pi K(o)}-g_\mu^{B}(o,\,o)\bigr]\e\cdot \kappa_4+o(\e)\Bigr),\end{equation}where $\kappa_4> 0$ is a universal constant that will be precised later.
From \cite[(37) and (38)]{Br1} we have\begin{eqnarray}
\int_{B}\phi_{\e}^6(x)\,\mathrm{d}x=\kappa_1\Bigl[1+8\pi g_\mu^{B}(o,\,o)\frac{\kappa_2}{\kappa_1}\e \Bigr]^3+o(\e)\nonumber\quad\\
=\kappa_1\Bigl[1+24\pi g_\mu^{B}(o,\,o)\frac{\kappa_2}{\kappa_1}\e \Bigr]+o(\e),
\label{B0}\,\,\,\\
\int_{B}|\nabla \phi_{\e}|^2-\mu \phi_{\e}^2=3\kappa_1\bigl[1+4\pi g_\mu^{B}(o,\,o)\cdot \frac{\kappa_2}{\kappa_1}\e\bigr] + o(\e),\label{B00}
\end{eqnarray}
where $\kappa_1=\int_{\mathbb{R}^3}
1/(1+|x|^2)^{3}\,\mathrm{d}\,x$ and
$\kappa_2=\int_{\mathbb{R}^3} 1/(1+|x|^2)^{\frac{5}{2}}\,\mathrm{d}\,x.$ Thus,
we are left with the
quantity $\int_B f_1(|x|)\phi_{\e}^6(x)\,\mathrm{d}\,x.$ Set
\begin{equation}\label{B1}h_{\e}:=\bigl(\phi_{\e}-U_{\e}\bigr)/\sqrt{\e},\end{equation}
and recall that we have from \cite[(39)]{Br1}
\begin{equation}\label{B2}h_{\e}\rightarrow 4\pi g_\mu^{B}(x,\,o)\qquad \mathrm{uniformly\,\,on}\,\,\bar{B}\,\,\mathrm{as\,\,}\e\rightarrow 0.\end{equation}
In particular, we get, for $\e_0> 0$ a constant small enough,
\begin{equation}\label{B3}\mathrm{sup}_{B}\mid h_{\e}\mid\leq M,\quad \forall\,\,0\leq\e \leq \e_0,\end{equation}
where $M> 0$ is a constant depending only on $\e_0.$ From other hand, by using the fact that $f_1 \in C^1\bigl([0,\,1]\bigr),$ we derive, for $\delta_1> 0$ a constant small enough, the existence of a constant $\delta > 0$ such that
\begin{equation}\label{B4}\mid f_1(t)-f_1^{'}(0)t\mid\leq \delta_1\, t,\quad \forall\,\,0\leq t\leq \delta.\end{equation}
Now, to estimate  $\int_B f_1(|x|)\phi_{\e}^6(x)\,\mathrm{d}\,x$ we argue as follows: The fact that $f_1(0)=0$ implies that
\begin{eqnarray}
\int_B f_1(|x|)\phi_{\e}^6(x)\,\mathrm{d}\,x=f_1^{'}(0)\underbrace{\int_B |x|\phi_{\e}^6(x)\,\mathrm{d}\,x}_{(II)}+\underbrace{\int_{B_{\delta}} \bigl(f_1(|x|)-f_1^{'}(0)|x|\bigr)\phi_{\e}^6(x)\,\mathrm{d}\,x}_{(I)}\label{B5}\\+\underbrace{\int_{B_{\delta}^c} \bigl(f_1(|x|)-f_1^{'}(0)|x|\bigr)\phi_{\e}^6(x)\,\mathrm{d}\,x}_{(III)},\qquad\,\,\,\quad\qquad\qquad\qquad\nonumber
\end{eqnarray}where $B_{\delta}:=B(o,\,\delta)$ and $B_{\delta}^c$ is its complementary in $\mathbb{R}^3.$ Combining \eqref{B1} and \eqref{B3}, and using the fact that $\sqrt{\e}= O(U_{\e})$ on $\bar{B},$ we obtain
\begin{eqnarray}
(III)=O\Bigl(\int_{B_{\delta}^c} U_{\e}^6(x)\,\mathrm{d}\,x\Bigr)\qquad\qquad\qquad\,\,\nonumber\\=\e^3\cdot \delta^{-3}\cdot O(1)=o(\e)\quad \mathrm{as\,\,}\e\rightarrow 0\label{B6},
\end{eqnarray}where $O(1)$ is a quantity independent of $\e$ and $\delta,$ and
\begin{eqnarray}
(II)=\int_B|x|U_{\e}^6(x)\,\mathrm{d}\,x+O\bigl(\,\int_B|x|U_{\e}^5(x)\sqrt{\e}\,\bigr)\,\mathrm{d}\,x\nonumber\\
=\int_{\mathbb{R}^3}|x|U_{\e}^6(x)\,\mathrm{d}\,x+o(\e)\nonumber\,\,\,\,\,\quad\qquad\qquad\qquad\\
=\kappa_3\cdot \e + o(\e),\quad\,\qquad\qquad\qquad\qquad\qquad\quad\label{B7}
\end{eqnarray}
where $\kappa_3:=\int_{\mathbb{R}^3}|x|/(1+|x|^3)^3\,\mathrm{d}\,x.$ Finally, from \eqref{B4} and the estimate of $(II)$ we get
\begin{eqnarray}
(I)\leq \delta_1\int_{B_{\delta}} |x|\phi_{\e}^6(x)\,\mathrm{d}\,x\Bigr)\nonumber\qquad\qquad\qquad\qquad\\=\delta_1\cdot \e \cdot O(1)=o(\e)\quad \mathrm{as\,\,}\delta_1\rightarrow 0,
\qquad\quad\label{B8}.
\end{eqnarray}where $O(1)$ is a quantity independent of $\e$ and $\delta_1.$
Combining \eqref{B5}--\eqref{B8}, we obtain
\begin{eqnarray}
\int_B f_1(|x|)\phi_{\e}^6(x)\,\mathrm{d}\,x=f_1^{'}(0)\cdot\kappa_3\cdot \e +o(\e).\label{B9}
\end{eqnarray}
Now, by combining \eqref{bahri}, \eqref{B0}, \eqref{B00} and \eqref{B9}, and using the fact that $S=3\cdot \kappa_1^{2/3}$ and $\kappa_2=(4/3)\kappa_3,$ we get, for $\e> 0$ small enough, $\phi_{\e}\in M_5$ and the following test
\begin{eqnarray}
J\bigl(\phi_{\e}\bigr)=\frac{3\kappa_1\bigl[1+4\pi g_\mu^{B}(o,\,o)\cdot \frac{\kappa_2}{\kappa_1}\e\bigr] + o(\e)}{\Bigl(K(o)\kappa_1\bigl[1+24\pi g_\mu^{B}(o,\,o)\frac{\kappa_2}{\kappa_1}\bigr]\cdot \e+\eta f_1^{'}(0)\kappa_3\cdot \e+o(\e)\Bigr)^{\frac{1}{3}}}\nonumber\\
=\frac{S}{\bigl(K(o)\bigr)^{\frac{1}{3}}}\Bigl(1-\bigl[4\pi g_\mu^B(o,\,o)+\frac{\eta f_1^{'}(0)}{4K(o)}\bigr]\frac{\e \cdot \kappa_2}{\kappa_1}+o(\e)\Bigr).\,\nonumber
\end{eqnarray}Thus, the claim \eqref{boucheche} follows with $\kappa_4=4\pi\cdot\kappa_2/\kappa_1.$ By combining \eqref{brezis} and \eqref{boucheche} and taking $\e> 0$ small enough, we obtain \eqref{lions} with $u=\phi_{\e}.$ This finishes the proof of the sufficiency.\\\\
{\it Necessity of the condition \eqref{brezis}:} Arguing by
contradiction, assuming that the problem $\mathbf{(BN)_{\eta}^1}$ has a
solution $u$ under the
condition\begin{equation}\label{condition}-\frac{f_1^{'}(0)\eta}{16\pi K(o)}\geq g_\mu^{B}(o,\,o).\end{equation}
By a result of Gidas--Ni--Nirenberg
\cite[Theorem 1$'$]{GNNir}, $\mathbf{(K_{\eta}^1)}$
implies that $u$ is necessarily spherically symmetric. We write
$u(x)=:u(t),$ where $t=|x|\in [0,\,1].$ To prove such kind of result, we need the Pohozaev-type identity \cite [(2.18)]{Bo1} which remains valid in three dimentional case: For any smooth function $\psi$ on $[0,\,1]$ with $\psi(0)=0,$ we have
\begin{align}\label{pohozaev}\begin{split}
&-\frac{1}{2}|u'(1)|^2\psi(1)+\int_0^1u^{2}\bigl[\mu\psi'(t)+\frac{1}{4}\psi^{(3)}(t)\bigr]t^2\,\mathrm{d}\,t\\&=
\int_0^1u^{6}\Bigl[-\frac{1}{6}tf'(t)\psi(t)+\frac{2}{3}f(t)\bigl(\psi(t)-t\psi'(t)\bigr)\Bigr]t\,\mathrm{d}\,t.
\end{split}\end{align}Now to get a contradiction, we need to look for a suitable function $\psi$ such that, for any $0\leq t\leq 1,$
\begin{eqnarray}
\psi(t)\geq 0,\qquad\qquad\qquad\,\,\,\quad\label{z1}\\
\mu\psi'(t)+\frac{1}{4}\psi^{(3)}(t)=0,\qquad\qquad\qquad\quad\,\,\,\label{z2}\\-f^{'}(t)t\psi(t)+4f(t)\bigl(\psi(t)-t\psi^{'}(t)\bigr)\geq 0,\qquad\qquad\qquad\quad\,\,\,\label{z3}\\-f^{'}(t_0)t_0\psi(t_0)+4f(t_0)\bigl(\psi(t_0)-t_0\psi^{'}(t_0)\bigr)> 0\quad \mathrm{for\,some\,\,t_0> 0}.\label{z3'}
\end{eqnarray}
To this end, we distinguish two cases:\\
{\bf Case 1}: $0< \mu \leq (\pi^2/4).$ In this case, we take, for any $0\leq t\leq 1,\,\,\psi(t)=\sin(2\sqrt{\mu}t)=:\psi_1(t).$
The claims \eqref{z1} and \eqref{z2} are clearly satisfied by the function $\psi_1(t).$ Combining $\mathbf{(K_{\eta}^1)}$ and \eqref{z1}, we get\begin{eqnarray}
-f^{'}(t)t\psi(t)+4f(t)\bigl(\psi(t)-t\psi^{'}(t)\bigr)\geq 4f(t)\bigl((\psi(t)-t\psi^{'}(t)\bigr).\label{z3''}
\end{eqnarray}
The claim \eqref{z3} follows from $\mathbf{(K_{\eta}^1)},$ \eqref{z3''} and the fact that \begin{equation}\label{z0}
\psi_1(t)-t\psi_1^{'}(t)\geq 0,\quad \forall\,\, 0\leq t\leq 1.
\end{equation}Finally, from \eqref{z3''} we derive that
\begin{eqnarray}
\lim_{t\rightarrow 0}\frac{-f^{'}(t)t\psi(t)+4f(t)\bigl(\psi(t)-t\psi^{'}(t)\bigr)}{t^3} \geq \lim_{t\rightarrow 0}\frac{4f(t)\bigl(\psi(t)-t\psi^{'}(t)\bigr)}{t^3}\nonumber\\
=\frac{32}{3}K(o)\mu\sqrt{\mu}> 0,\qquad\quad\nonumber
\end{eqnarray}
and then the claim \eqref{z3'} follows.\\\\
{\bf Case 2}: $(\pi^2/4)< \mu < \pi^2.$ Set, for any $0\leq t\leq 1,\,\,\psi_2(t):=1-\cos(2\sqrt{\mu}t)$ and define$$\psi(t)=\psi_1(t)-\frac{\cos(\sqrt{\mu})}{\sin(\sqrt{\mu})}\psi_2(t),\quad \forall\,\,0\leq t\leq1.$$
A straightforward calculation shows that the function $\psi(t)$ satisfies the claims \eqref{z1} and \eqref{z2}. To check the claim \eqref{z3} we need to discuss two regions of the interval $[0,\,1].$ To this end, observe that we have from the condition $\mathbf{(K)}$ the existence of two constants $\delta_0> 0$ and $M_0\geq 0$ such that, for any $0\leq t\leq \delta_0,$
\begin{equation}\label{zakaria4}f_1^{'}(t)\leq M_0\cdot t+f_1^{'}(0).
\end{equation}
Next by studying the variations of the function $\bigl(\psi(t)-t\psi^{'}(t)\bigr),$ we derive the existence of a unique $t_{\mu}$ in $]0,\,1/2]$ such that
\begin{eqnarray}
\psi(t_{\mu})-t_{\mu}\psi^{'}(t_{\mu})=0,\qquad\qquad\qquad\,\,\,\,\,\nonumber\\
\psi(t)-t\psi^{'}(t)\geq 0,\quad \forall\,\,t_{\mu}\leq t\leq 1,\label{z5}\\
\psi(t)-t\psi^{'}(t)\leq0 ,\quad \forall\,\,0\leq t\leq t_{\mu}.\,\label{z6}
\end{eqnarray}
Observe that $t_{\mu}\rightarrow 0$ as $\bigl(\eta/ K(o)\bigr)\rightarrow 0,$ then we can choose $\eta_1> 0$ a constant such that, for any $0\leq \eta \leq \eta_1,$ we have $t_\mu\leq \delta_0.$ By combining $\mathbf{(K_{\eta}^1)},$ \eqref{z1} and \eqref{z5}, the claim \eqref{z3} follows for any $t_{\mu}\leq t\leq1.$ Now, we are left with the region $[0,\,t_{\mu}].$ In this case, we develop the left hand side of \eqref{z3} as follows:
\begin{eqnarray}
-f^{'}(t)t\psi(t)+f(t)\cdot\bigl(\psi(t)-t\psi^{'}(t)\bigr)=
\underbrace{4\bigl(\psi_1(t)-t\psi_1^{'}(t)\bigr)f(t)+\eta\bigl(f_1^{'}(0)-f_1^{'}(t)\bigr)t\psi_1(t)}_{(IV)}\,\,\,\,\nonumber\\
\underbrace{-\frac{\cos(\sqrt{\mu})}{\sin(\sqrt{\mu})}\Bigl[4\bigl(\psi_2(t)-t\psi_2^{'}(t)\bigr)\eta f_1(t)+f^{'}(t)t\psi_2(t)\Bigr]}_{(V)}\label{z000}\\ \underbrace{-f_1^{'}(0)\eta t\psi_1(t)-4K(o)\frac{\cos(\sqrt{\mu})}{\sin(\sqrt{\mu})}\bigl(\psi_2(t)-t\psi_2^{'}(t)\bigr)}_{(VI)}.\,\nonumber
\end{eqnarray}
By using the fact that $0\leq t_{\mu}\leq 1/2,$ we derive from \eqref{z0} and \eqref{z6} that
\begin{equation*}\label{z7}\psi_2(t)-t\psi_2^{'}(t)\leq 0,\quad \forall\,\, 0\leq t\leq t_{\mu}.
\end{equation*}
This, together with $\mathbf{(K_{\eta}^1)}$ and the fact that $\pi^2/4< \mu < \pi^2,$ implies that
\begin{equation}\label{z8}(V)\geq 0.
\end{equation}Next we derive from \eqref{z00} and \eqref{condition} that
\begin{eqnarray}
(VI)\geq
-4K(o)\frac{\cos(\sqrt{\mu})}{\sin(\sqrt{\mu})}\bigl(\sqrt{\mu}\,t\psi_1(t)+\psi_2(t)-t\psi_2^{'}(t)\bigr)\nonumber\,\,\,\,\,\,\,\,\quad\qquad\qquad\\
\geq-8K(o)\frac{\cos(\sqrt{\mu})}{\sin(\sqrt{\mu})}\psi_1(\frac{t}{2})\bigl(\psi_1(\frac{t}{2})-\frac{t}{2}\psi_1^{'}(\frac{t}{2})\bigr)\label{z9'}\,\,\,\,\quad\qquad\qquad\quad\quad\,\,\\
\geq 0\quad \mathrm{from\,\,the\,\,study\,\,of\,\,{\bf case \,1}\,\,and\,the\,fact\,that\,\,}\frac{\pi^2}{4}< \mu< \pi^2\label {z9}.
\end{eqnarray}
Finally to deal with $(IV)$ we use \eqref{zakaria4} which implies that
\begin{equation}\label{z10} (IV)\geq 4\bigl(\psi_1(t)-t\psi_1^{'}(t)\bigr)f(t)-\eta M_0t^2\psi_1(t).
\end{equation}
Observe that since $0\leq t\leq 1/2,$ then by reason of continuity there exist two constants $\kappa_5,\,\kappa_6> 0$ independents of $\mu$ such that
\begin{equation}\label{z11}  t^3\cdot \kappa_5\leq \bigl(\psi_1(t)-t\psi_1^{'}(t)\bigr),\,t^2\psi_1(t)\leq t^3\cdot \kappa_6.
\end{equation}
Let $\eta_2> 0$ be a constant such that, for any $0\leq \eta \leq \eta_2,$ we have\begin{eqnarray}
K(o)+\eta_2f_1(\delta_0)=:\kappa_7> 0,\label{z12}\\
4\kappa_5\kappa_7-\eta M_0\kappa_6 \geq 0.\,\label{z13}
\end{eqnarray}By combining $\mathbf{(K_{\eta}^1)}$ and \eqref{z10}--\eqref{z13}, we derive that
\begin{equation}\label{z14}(IV)\geq 0.
\end{equation}
Thus, for any $0\leq \eta \leq \bar{\eta}:=\mathrm{min}(\eta_1,\,\eta_2),$ the claim \eqref{z3} follows from \eqref{z000}, \eqref{z8}, \eqref{z9} and \eqref{z14}. Finally, by combining \eqref{z000}--\eqref{z9'} and \eqref{z14}, we obtain \begin{eqnarray}
-f^{'}(t)t\psi(t)+4f(t)\bigl(\psi(t)-t\psi^{'}(t)\bigr)\geq -8K(o)\frac{\cos(\sqrt{\mu})}{\sin(\sqrt{\mu})}\psi_1(\frac{t}{2})\bigl(\psi_1(\frac{t}{2})-\frac{t}{2}\psi_1^{'}(\frac{t}{2})\bigr).\nonumber
\end{eqnarray}This, together with the fact that $$\lim_{t\rightarrow 0}\frac{\psi_1(\frac{t}{2})\bigl(\psi_1(\frac{t}{2})-\frac{t}{2}\psi_1(\frac{t}{2})\bigr)}{t^4} =\frac{\mu^2}{3}>0,$$ implies the claim \eqref{z3'}.
 As a conclusion, the claims \eqref{z1}--\eqref{z3'} follow for any $0< \mu < \pi^2$ and any $0\leq \eta \leq \bar{\eta}$ satisfying \eqref{condition}. This contradicts \eqref{pohozaev}, and then the proof of the necessity of the condition \eqref{brezis} follows for any
$0\leq \eta \leq \bar{\eta}.$ This finishes the proof of Theorem \ref{optimal}.

\end{document}